\documentclass[a4paper,12pt]{article}

\usepackage[margin=2.5cm,bmargin=3cm,tmargin=2.8cm]{geometry}
\usepackage{amsmath}
\usepackage{amssymb}
\usepackage{tensor}
\usepackage{latexsym}
\usepackage{enumerate}
\usepackage{bbm}
\usepackage{amsthm}
\usepackage{cite}

\newcommand{\ima}{\mathbbmtt{i}}

\title{\Large{\textbf{Orthogonality of the Ferrers' Associated Legendre Functions of the Second Kind with Imaginary Argument}}}
\author{N. Dimakis\thanks{nsdimakis@gmail.com} \\ \small{Instituto de Ciencias F\'isicas y Matem\'aticas}\\ \small{Universidad Austral de Chile, 5090000 Valdivia, Chile}}
\date{}

\begin{document}
\maketitle
\numberwithin{equation}{section}

\begin{abstract}
  In this work we study the associated Legendre functions of the second kind with a purely imaginary argument $Q^k_\ell(\ima\, x)$. We derive the conditions under which they provide a set of square integrable functions when $x \in \mathbb{R}$ and we prove the relevant orthogonality relation that they satisfy.
\end{abstract}

\section{Introduction}

The importance of the associated Legendre equation is well known in applied mathematics. Its significance is also evident in mathematical physics as it appears in various different areas of study; from geodesy \cite{Heisk} to quantum mechanics \cite{Edmonds}. In the latter it usually emerges through a separation of values process over the eigenvalue equation for the Casimir invariant of the $\mathfrak{so}(3)$ algebra of the angular momentum operators \cite{Edmonds}, that leads to
\begin{equation} \label{Leqeq}
  \frac{d}{d z} \left[ (1-z^2) \frac{dY(z)}{dz}\right] + \left[ \ell(\ell+1) - \frac{k^2}{1-z^2}\right] Y(z)=0.
\end{equation}
The solution of \eqref{Leqeq} is expressed as a linear combination of the two associate Legendre functions of the first $P^k_\ell(z)$ and second kind $Q^k_\ell(z)$,
\begin{equation}
  Y(z) = c_1\, P^k_\ell(z) + c_2\, Q^k_\ell(z),
\end{equation}
which can be defined by means of a more general function:
\begin{subequations} \label{LegPQdef}
  \begin{align}
    P^k_{\ell}(z) =& \frac{(1+z)^{k/2}}{(1-z)^{k/2}} {}_2 \tilde{F}_1 (-\ell, \ell+1;1-k ; \frac{1-z}{2}) \\
    \begin{split}
      Q^k_{\ell}(z) =& \frac{\pi}{2 \sin(k \pi)} \Bigg[ \frac{(1+z)^{k/2}}{(1-z)^{k/2}} \cos(k\pi) {}_2 \tilde{F}_1 (-\ell, \ell+1;1-k ; \frac{1-z}{2})  \\
      & - \frac{\Gamma(k+\ell+1)}{\Gamma(-k+\ell+1) }\frac{(1-z)^{k/2}}{(1+z)^{k/2}} {}_2 \tilde{F}_1 (-\ell, \ell+1;k+1 ; \frac{1-z}{2}) \Bigg]
    \end{split}
  \end{align}
  \end{subequations}
where $\Gamma(\alpha)$ represents the Gamma function, while
\begin{equation} \label{normF}
  {}_2 \tilde{F}_1 (a, b;c ; z) = \sum_{s=0}^{+\infty} \frac{(a)_s (b)_s}{\Gamma (c+s) s!} z^s, \quad |z|<1
\end{equation}
is the regularized Gauss hypergeometric function, in which also appears the Pochhammer symbol $(x)_n = \Gamma(x+n)/\Gamma(x)$.

At this point we have to note that $P^k_{\ell}(z)$ and $Q^k_{\ell}(z)$ as defined by \eqref{LegPQdef} are the analytic continuations of Ferrers' functions \cite{NIST} and are different from those obtained by the definition of the associated Legendre functions by Hobson \cite{Hobson}. Wherever needed we shall denote the latter pair of functions as $\mathcal{P}^k_{\ell}(z)$ and $\mathcal{Q}^k_{\ell}(z)$.

There exists an extensive bibliography regarding the properties of \eqref{LegPQdef} \cite{NIST,Hobson,Abr,Magn,FOlver}. Most of those that are well known refer to the associate Legendre function of the first kind $P^k_\ell(z)$, with the most recognizable being of course the orthogonality relation
\begin{equation}
  \int_{-1}^{1}\!\! P^k_m(x) P^k_n(x) dx = \frac{2(k+m)! \delta_{mn}}{(2m+1)(m-k)!}, \quad k,m,n \in \mathbb{N}.
\end{equation}
Other orthogonality relations regarding $P^k_\ell(z)$ in terms of Dirac's delta function when the order or the degree is complex are also known, see for example \cite{Nostrand,Gotze,Kalmykov,Bielski}. In this work however, we concentrate on $Q^k_\ell(z)$ and to the orthogonality relation to which it leads when the argument is purely imaginary. In particular, we prove that when $z=\ima\, x$, $x \in \mathbb{R}$, the following relation holds
\begin{equation} \label{orthrel}
  \int_{-\infty}^{+\infty}\!\! Q^k_\ell (\ima\, x)Q^k_{\ell'} (\ima\, x)dx =
  \begin{cases}
    (2\ell)! \pi \left[ \prod_{i=\ell+2}^{k}(\ell+i)(\ell-i+1) \right] \delta_{\ell\ell'} , & k>\ell+1 \\ \\
    (2\ell)! \pi \delta_{\ell\ell'} , & k=\ell+1
  \end{cases}
\end{equation}
whenever $k\in \mathbb{Z}_+$, $\ell \in \mathbb{N}$, with $k>\ell$. The definite integral of the left hand side denotes the improper Riemann integral defined as
\begin{equation}
\int_{-\infty}^{+\infty}\!\! F(x) dx = \underset{a\rightarrow +\infty}{\lim} \int_{-a}^{a}\!\! F(x)dx.
\end{equation}

The consideration of an imaginary argument for the independent variable $z= \ima \, x$, turns \eqref{Leqeq} into
\begin{equation} \label{Legeqx}
  \frac{d}{d x} \left[ (1+x^2) \frac{d Y(x)}{dx}\right] - \left[ \ell(\ell+1) - \frac{k^2}{1+x^2}\right] Y(x)=0
\end{equation}
with the solution being expressed of course as any linear combination of $P^k_{\ell}(\ima\, x)$ and $Q^k_{\ell}(\ima\, x)$. This equation is our starting point in the consecutive analysis. The structure of the paper is the following: In section 2 we briefly summarize some of the properties that are to be used. In particular, we are interested in the behaviour of $Q^k_\ell (\ima\, x)$ at the boundary $x\rightarrow \pm \infty$ and for which values the function vanishes. In section 3 we start from \eqref{Legeqx} and exhibit under which conditions $Q^k_{\ell}(\ima\, x)$ is orthogonal in $x \in (-\infty,+\infty)$ for different values of $\ell$, while in the final section we prove that it is square integrable and derive expression \eqref{orthrel}.

\section{Some relevant properties of $Q^k_{\ell}(z)$}

First of all, let us investigate how $Q^k_{\ell}(z)$, with $z=\ima\, x$, behaves as a function in relation to its order and degree. In what follows we assume that the order of the function is a positive integer i.e. $k \in \mathbb{Z}_{+}$ and we distinguish the following two cases:
\begin{enumerate}[(1)]
  \item \label{caseQ1} $\ell \in \mathbb{N}$\\
  Then, $Q^k_{\ell}(z)$ can be expressed with the help of the subsequent relations
  \begin{subequations} \label{LegPQ}
    \begin{align} \label{intdefQkl}
      Q^k_{\ell}(z) & = (-1)^k(1-z^2)^{k/2} \frac{d^k}{dz^k}Q_{\ell}(z) \\ \label{intdefQsolol}
      Q_{\ell}(z) & = \frac{1}{2} P_{\ell}(z) \ln\left(\frac{1+z}{1-z}\right) - \sum_{j=1}^{\ell}\frac{P_{j-1}(z)P_{\ell-j}(z)}{j} \\
      P_{\ell}(z) & = \frac{1}{2^{\ell} \ell !} \frac{d^\ell}{dz^\ell}\left[\left(z^2-1\right)^\ell\right],
    \end{align}
  \end{subequations}
  where $P_{\ell}$, $Q_{\ell}$ are the Legendre functions of the first and second kind respectively.

  It is true that the associated Legendre function of the second kind $\mathcal{Q}_\ell^k(z)$, as defined by Hobson, is more frequently encountered in the literature. In that regard, one can encounter the following interesting set of properties \cite{Susc}:
  \begin{equation} \label{Qpol}
    \mathcal{Q}_\ell^k(z) = \begin{cases}
                    \frac{(z^2-1)^{-k/2}}{(\ell-k)!}\left[\frac{1}{2} A^{(1)}_{k\ell}(z) \ln\left(\frac{z+1}{z-1}\right)+ A^{(2)}_{k\ell}(z)\right], & \ell \geq k \\ \\
                    \frac{(-1)^k}{(z^2-1)^{k/2}}A^{(2)}_{k\ell}(z), & \ell < k
                  \end{cases}
  \end{equation}
  where $A^{(i)}_{k\ell}(z)$, $i=1,2$, are polynomials in $z$, with the following properties:
  \begin{equation} \label{propA}
    \begin{split}
       &   A^{(1)}_{k\ell}(z)  \begin{cases}
                          \text{is of order } \ell+k, & \; \text{when } \ell \geq k \\
                          =0, & \; \text{when } \ell < k
                       \end{cases}
  \\
       &   A^{(2)}_{k\ell}(z)  \begin{cases}
                          \text{is of order } \ell+k-1, & \; \text{when } \ell \geq k \\
                          \text{is of order } k-\ell-1, & \; \text{when } \ell < k
                       \end{cases}
    \end{split}
  \end{equation}
  and additionally
  \begin{itemize}
    \item The coefficients of the powers of $z$ in $A^{(i)}(z)$ are all integers.
    \item If the order of the polynomial is even (odd) then all the powers of $z$ are even or odd.
  \end{itemize}
  Hence, for $k>\ell$, it follows that when $k-\ell-1$ is even, $\mathcal{Q}^k_\ell(\ima\, x)$ is a real function, else it is purely imaginary. What is more, as can be seen from \eqref{Qpol} and \eqref{propA}, $\mathcal{Q}^k_\ell(\ima\, x)$ and all its derivatives are finite for every $x \in \mathbb{R}$.

  It so happens that for the range of $k$ and $\ell$ considered here, the difference between $Q^k_\ell(z)$ and $\mathcal{Q}^k_\ell(z)$ reduces to that of a multiplicative constant. The latter is of the form $e^{n \ima \pi/2}$, $n\in \mathbb{N}$ and depends on a combination of the order $k$ and the sign of $z$ if $z\in \mathbb{R}$ or the sign of $\mathrm{Im}(z)$ if $z$ is complex. In particular it holds that \cite{Mathematica}
  \begin{equation} \label{relationbttwo}
    Q^m_n (z) = e^{\pm \ima \pi m/2} \mathcal{Q}^m_n(z), \quad m, n \in \mathbb{N}; m>n; \mathrm{Im} z \gtrless 0
  \end{equation}
  where $[s]$ denotes the larger integer that is equal or less than $s$. As a result, we can also utilize these properties to characterize the general behaviour of $Q^k_{\ell}(\ima\, x)$. Hence, we can state that, when the inequality $k>\ell$ holds, $Q^k_{\ell}(\ima\, x)$ is a solution of \eqref{Legeqx} which vanishes at infinity and at the same time does not posses any singularities for $x \in \mathbb{R}$. Although this is neither a sufficient nor an adequate condition for a function to be square integrable, it is a first indication of a behaviour that may lead to such a result.

  \item $\ell \in \mathbb{R}$

  In this case let us use the following expression that is valid for values of $z$ outside the unit circle
  \begin{equation} \label{Qserexp}
    \begin{split}
       Q^k_\ell (z) = & \frac{2^{-(\ell+2)}e^{\ima k \pi}\sqrt{\pi}(1-z^2)^{k/2}}{\cos(\ell \pi) z^{k+\ell+1}}  \\
       &\times \left[ \ima\, 2^{2\ell+1} \Gamma(k-\ell) \sin\left((k-\ell)\pi\right) z^{2\ell+1} \sum_{s=0}^{+\infty}\frac{\left(\frac{k-\ell}{2}\right)_s \left(\frac{k-\ell-1}{2}\right)_s}{\Gamma(s-\ell+\frac{1}{2})s! z^{2 s}} \right.  \\
          & \left. + \left( \cos\left((k+\ell)\pi\right)+e^{\ima(\ell-k)\pi}\right)\Gamma(k+\ell+1) \sum_{s=0}^{+\infty}\frac{\left(\frac{k+\ell+1}{2}\right)_s \left(\frac{k+\ell+2}{2}\right)_s}{\Gamma(s+\ell+\frac{3}{2})s! z^{2 s}} \right].
    \end{split}
  \end{equation}
  With the help of \eqref{Qserexp} let us examine what happens to $Q^k_\ell (\ima\, x)$ as $x\rightarrow \pm \infty$ and in particular possible values of $\ell$ for which the function vanishes at the boundary. Given the fact that  only the $s=0$ terms are important in the sums of \eqref{Qserexp} as $z$ tends to infinity, we can distinguish the following two cases:
  \begin{enumerate}
   \item Firstly, $\ell \geq 0$, then in order to have a zero at infinity we must demand $k-\ell \in \mathbb{Z}_{+}$, which is again what we got in the previous investigation.

   \item On the other hand, when $\ell <0$, we must necessarily enforce $\ell>-1$.
  \end{enumerate}

  As a result we can state that $Q^k_\ell (\ima\, x)$ vanishes at the boundary $\pm\infty$ either if $k>\ell \in \mathbb{N}$ or when $\ell \in (-1,0)$.
\end{enumerate}

\section{Orthogonality}

We start from equation \eqref{Legeqx} and consider the two following relations that hold identically for $Q^k_\ell(\ima\, x)$ and $Q^m_n(\ima\, x)$
\begin{subequations}\label{ortheqs}
  \begin{align}\label{ortheq1}
    & \frac{d}{dx}\left[(1+x^2)\frac{d Q^k_\ell}{dx}\right] -\left[\ell(\ell+1)-\frac{k^2}{1+x^2}\right]Q^k_\ell =0 \\ \label{ortheq2}
    & \frac{d}{dx}\left[(1+x^2)\frac{d Q^m_n}{dx}\right] -\left[n(n+1)-\frac{m^2}{1+x^2}\right]Q^m_n =0 .
  \end{align}
\end{subequations}
If we multiply the first and the second with $Q^m_n(\ima\, x)$, $Q^k_\ell(\ima\, x)$ respectively and succeedingly subtract the one from the other we arrive at
\begin{equation*}
  \left[(\ell-n)(\ell+n+1)-\frac{k^2-m^2}{1+x^2}\right]Q^k_\ell Q^m_n = \frac{d}{dx}\left[(1+x^2)\left(Q^m_n\frac{d Q^k_\ell}{dx}- Q^k_\ell \frac{d Q^m_n}{dx}\right)\right],
\end{equation*}
which, for $k=m$ and $\ell\neq n$ becomes
\begin{equation} \label{orthcond}
  (\ell+n+1)\int_{-\infty}^{+\infty}\!\! Q^k_\ell(\ima\, x) Q^k_n(\ima\, x) dx = \left[\frac{1+x^2}{\ell-n}\left(Q^k_n\frac{d Q^k_\ell}{dx}- Q^k_\ell \frac{d Q^k_n}{dx}\right)\right]_{-\infty}^{+\infty}.
\end{equation}
If we want an orthogonality condition to be satisfied when $\ell\neq n$ the boundary term that appears on the right hand of \eqref{orthcond} must be zero. Until now we have deduced the two subsequent possibilities from the previous section:
\begin{enumerate}
  \item $k\in \mathbb{Z}$, $\ell \in (-1,0)$
  \item $k\in \mathbb{Z}_{+}$, $\ell \in \mathbb{N}$ and $k>\ell$.
\end{enumerate}
For the first case let us check the behaviour of the product $x^2\, Q^k_\ell \, \frac{d Q^m_n}{dx}$ for the leading terms as $x\rightarrow \pm \infty$, which can be deduced with the help of \eqref{Qserexp}:
\begin{equation*}
  \begin{split}
    x^2 Q^k_\ell \frac{d Q^k_n}{dx} & \sim x^2\left(\frac{\alpha_1}{x^{|\ell|}}+\frac{\alpha_2}{x^{1-|\ell|}}\right) \left(\frac{\beta_1}{x^{1+|n|}}+\frac{\alpha_2}{x^{2-|n|}}\right) \\
    & \sim \frac{\gamma_1}{x^{|n|+|\ell|-1}} + \frac{\gamma_2}{x^{|n|-|\ell|}} + \frac{\gamma_3}{x^{|\ell|-|n|}} + \frac{\gamma_1}{x^{1-|n|-|\ell|}},
  \end{split}
\end{equation*}
where the $\alpha_i$, $\beta_i$ and $\gamma_i$'s are constants. It is clear that the expression diverges at the boundary for all $n,\ell \in (-1,0)$.

On the contrary, when we turn to the last case and remember \eqref{Qpol} together with properties \eqref{propA} we see that $Q^k_\ell(\ima\, x)$ decays as $x^{-\ell-1}$, thus
\begin{equation*}
  x^2\, Q^k_\ell \, \frac{d Q^k_n}{dx} \sim x^2 \frac{1}{x^{\ell+1}}\frac{1}{x^{n+2}} =\frac{1}{x^{\ell+n+1}}
\end{equation*}
which vanishes on the boundary, since now $\ell,n \geq 0$. Hence we can see that, functions $Q^k_\ell(\ima\, x)$, for $k \in \mathbb{Z}_{+}$, $\ell \in \mathbb{N}$, with $k>\ell$ form an orthogonal set in the region $x \in (-\infty,+\infty)$.

\section{Square integrability}

By having proven that
\begin{equation} \label{orthzero}
  \int_{-\infty}^{+\infty}\!\! Q^k_\ell(\ima\, x) Q^{k}_{n}(\ima\, x) = 0,\quad \ell \neq n
\end{equation}
for the given range of values of $k,\ell$ and $n$, we need only investigate what happens when $\ell = n$, case where equation \eqref{orthcond} is not applicable. We work in a similar way to the derivation of the normalization of the Legendre polynomials of the first kind $P^k_\ell(x)$ when $x\in[-1,1]$ \cite{Tang}. At first we use definition \eqref{intdefQkl} to write
\begin{equation*}
  Q^k_\ell (\ima\, x)Q^k_\ell (\ima\, x) = \frac{(1+x^2)^k}{\ima^{2k}}\frac{d^k Q_\ell(\ima\, x)}{dx^k} \frac{d^k Q_\ell(\ima\, x)}{dx^k}.
\end{equation*}
If we integrate the previous relation, we get
\begin{equation} \label{orthint1}
  \begin{split}
    \int_{-\infty}^{+\infty}\!\! Q^k_\ell (\ima\, x)Q^k_\ell (\ima\, x)dx  = & \frac{1}{\ima^{2k}}\left[(1+x^2)^k \frac{d^k Q_\ell(\ima\, x) }{dx^k}\frac{d^{k-1}Q_\ell(\ima\, x)}{dx^{k-1}}\right]_{-\infty}^{+\infty} \\
    & - \frac{1}{\ima^{2k}}\int_{-\infty}^{+\infty}\!\! \frac{d}{dx}\left((1+x^2)^k \frac{d^k Q_\ell(\ima\, x)}{dx^k}\right)\frac{d^{k-1} Q_\ell(\ima\, x)}{dx^{k-1}} dx,
  \end{split}
\end{equation}
where the surface term of the right hand side of \eqref{orthint1} is zero, since the expression in the brackets decayes at infinity like $x^{-2\ell-1}$. Let us now see what happens with the second term. At first we need to calculate
\begin{equation} \label{orthint2}
  \frac{d}{dx}\left((1+x^2)^k \frac{d^k}{dx^k}Q_\ell(\ima\, x)\right) = (1+x^2)^k \frac{d^{k+1}}{dx^{k+1}}Q_\ell(\ima\, x) + 2k x (1+x^2)^{k-1} \frac{d^k}{dx^k}Q_\ell(\ima\, x),
\end{equation}
then we need to remember that the Legendre function of the second kind $Q_\ell(\ima\, x)$ satisfies the differential equation
\begin{equation}\label{LegeqQ}
  (1+x^2)\frac{d^2}{dx^2}Q_\ell(\ima\, x)+2 x\frac{d}{dx}Q_\ell(\ima\, x)-\ell(\ell+1)Q_\ell(\ima\, x)=0.
\end{equation}
By differentiating $(k-1)$ times with respect to $x$ and with the help of Leibnitz's formula
\begin{equation}\label{Leib}
  \frac{d^k}{dx^k}\left[A(x)B(x)\right] = \sum_{m=0}^{k}\frac{k!}{m!(k-m)!}\frac{d^{k-m} A(x)}{dx^{k-m}}\frac{d^m B(x)}{dx^m},
\end{equation}
we are led to the following relation
\begin{equation}\label{orthint3}
  (1+x^2)\frac{d^{k+1}}{dx^{k+1}}Q_\ell(\ima\, x)+2k x \frac{d^k}{dx^k}Q_\ell(\ima\, x) = (\ell+k)(\ell-k+1)(1+x^2)^{k-1}\frac{d^{k-1}}{dx^{k-1}}Q_\ell(\ima\, x).
\end{equation}
The left hand side of \eqref{orthint3} is equal to the right hand side of \eqref{orthint2}, hence we can deduce that
\begin{equation*}
  \begin{split}
    \int_{-\infty}^{+\infty}\!\! Q^k_\ell (\ima\, x)Q^k_\ell (\ima\, x)dx  = & \frac{(\ell+k)(\ell-k+1)}{-\ima^{2k}}
     \int_{-\infty}^{+\infty}\!\! (1+x^2)^{k-1}\frac{d^{k-1}Q_\ell(\ima\, x)}{dx^{k-1}}\frac{d^{k-1}Q_\ell(\ima\, x)}{dx^{k-1}}dx \\
     = & (\ell+k)(\ell-k+1) \int_{-\infty}^{+\infty}Q^{k-1}_\ell (\ima\, x)Q^{k-1}_\ell (\ima\, x)dx .
  \end{split}
\end{equation*}
By repeating the same procedure $k-\ell+1$ times we get
\begin{equation}\label{orthint4}
  \int_{-\infty}^{+\infty}\!\! Q^k_\ell (\ima\, x)Q^k_\ell (\ima\, x)dx = \left[\prod_{i=\ell+2}^{k}(\ell+i)(\ell-i+1) \right] \int_{-\infty}^{+\infty}\!\! Q^{\ell+1}_\ell (\ima\, x)Q^{\ell+1}_\ell (\ima\, x) dx.
\end{equation}
At this point we need only prove that the integral on the right hand side is bounded. Let us recall properties \eqref{propA}; for $k=\ell+1$, the polynomial $A^{(2)}_{k\ell}(z)$ is of zero-th order, hence \eqref{Qpol} implies
\begin{equation}\label{orthint5}
  Q^{\ell+1}_\ell(\ima\, x) = \frac{A_0(\ell)}{(1+x^2)^{(\ell+1)/2}},
\end{equation}
where $A_0(\ell)$ is defined by \eqref{orthint5} and its value with respect of $\ell$ given by \eqref{valueofAfinal}. Additionally, it is easy to verify that
\begin{equation}\label{orthint6}
  \int\!\! \frac{1}{(1+x^2)^{\ell+1}}dx = x\, {}_2F_1(\frac{1}{2},\ell+1;\frac{3}{2};-x^2) + \mathrm{const.}
\end{equation}
With the use of
\begin{equation} \label{regFrel}
  \begin{split}
    \frac{\sin\left((b-a)\pi\right)}{\pi}{}_2 \tilde{F}_1 (a,b;c;z) = \frac{(-z)^{-a}}{\Gamma(b)\Gamma(c-a)}{}_2 \tilde{F}_1(a,a-c+1;a-b+1;\frac{1}{z}) \\
    - \frac{(-z)^{-b}}{\Gamma(a)\Gamma(c-b)}{}_2 \tilde{F}_1(b,b-c+1;b-a+1;\frac{1}{z})
  \end{split}
\end{equation}
which holds outside the unit circle $|z|=1$ \cite{FOlver}, definition \eqref{normF} and the relation that connects the regularized ${}_2\tilde{F}_1$ with the ordinary Gauss hypergeometric function, i.e. ${}_2F_1(a,b;c;z)=\Gamma(z)\,  {}_2\tilde{F}_1(a,b;c;z)$, we deduce that
\begin{equation} \label{approx}
  {}_2F_1(\frac{1}{2},\ell+1;\frac{3}{2};-x^2)|_{x\rightarrow \pm \infty} = \frac{(-1)^\ell\, \pi}{\Gamma(\frac{1}{2}-\ell)}\left[ \frac{ \pi^{1/2}}{2\, \Gamma(\ell+1)|x|} - \frac{(-x^2)^{-(\ell+1)}}{2\, \Gamma(\frac{3}{2}+\ell)}\right] + O[\frac{1}{x^2}].
\end{equation}
Given that the smallest value of $\ell$ is zero - so that the second term in the previous bracket can be neglected as $x\rightarrow \pm \infty$ - from relations \eqref{orthint5}-\eqref{approx} we get
\begin{equation}\label{prefinalrel}
  \begin{split}
    \int_{-\infty}^{+\infty} \left(Q^{\ell+1}_\ell (\ima \, x)\right)^2 dx  & = A_0^2(\ell) \left[x\,\frac{(-1)^\ell \pi^{3/2}}{2\, \Gamma(\frac{1}{2}-\ell)\Gamma(\ell+1)|x|}\right]_{-\infty}^{+\infty} \\
    & = \frac{(-1)^\ell \pi^{3/2} A_0^2(\ell)}{\Gamma(\frac{1}{2}-\ell)\Gamma(\ell+1)}
  \end{split}
\end{equation}
which leads \eqref{orthint4} to be written as
\begin{equation} \label{finalrel}
  \int_{-\infty}^{+\infty}\!\! Q^k_\ell (\ima\, x)Q^k_\ell (\ima\, x)dx = \frac{(-1)^\ell\, \pi^{3/2}A_0^2(\ell)}{\Gamma(\frac{1}{2}-\ell)\, \ell !} \prod_{i=\ell+2}^{k}(\ell+i)(\ell-i+1),
\end{equation}
which as we see is a bounded expression. In the appendix that follows we explicitly calculate the value of $A_0(\ell)$ and prove it to be $A_0(\ell)=(-1)^{\ell+1} 2^\ell \ell!$. Thus, by its substitution in \eqref{finalrel} together with $\Gamma(\frac{1}{2}-\ell) = \frac{(-1)^\ell 2^{2\ell}\ell!}{(2\ell)!} \sqrt{\pi}$, we finally arrive at
\begin{equation} \label{finalrel2}
  \int_{-\infty}^{+\infty}\!\! Q^k_\ell (\ima\, x)Q^k_\ell (\ima\, x)dx = (2\ell)! \pi \prod_{i=\ell+2}^{k}(\ell+i)(\ell-i+1).
\end{equation}
As a result we can state that, the $Q^k_\ell (\ima\, x)$'s for $k\in \mathbb{Z}_{+}$, $\ell \in \mathbb{N}$ with $k>\ell$ form an orthogonal set of square integrable functions in the region $x\in (-\infty,+\infty)$. What is more, by combining \eqref{orthzero}, \eqref{prefinalrel} and \eqref{finalrel2} we have proven that \eqref{orthrel} holds.

\appendix

\section{Calculation of $A_0(\ell)$}

Let us proceed with the calculation of $A_0(\ell)$ by using definitions \eqref{LegPQ}. From the latter we can derive the relation
\begin{equation} \label{app1}
  Q^{\ell+1}_\ell (z) = \frac{(-1)^{\ell+1}}{2} (1-z^2)^{(\ell+1)/2} \frac{d^{\ell+1}}{dz^{\ell+1}}\left[P_\ell(z) \ln\left(\frac{1+z}{1-z}\right)\right],
\end{equation}
since the second term of \eqref{intdefQsolol} is a polynomial of order $\ell-1$ and its $\ell+1$ order derivative is zero. By application of Leibnitz's formula \eqref{Leib} we can write
\begin{equation} \label{app2}
  \frac{d^{\ell+1}}{dz^{\ell+1}}\left[P_\ell(z) \ln\left(\frac{1+z}{1-z}\right)\right] = \sum_{m=1}^{\ell+1} \frac{(\ell+1)!}{m! (\ell+1-m)!} \frac{d^{\ell+1-m}P_{\ell}(z)}{dz^{\ell+1-m}} \frac{d^m}{dz^m}\ln\left(\frac{1+z}{1-z}\right),
\end{equation}
where the summation starts from $m=1$ and not $m=0$ because $P_{\ell}(z)$ is a polynomial of order $\ell$ and its $\ell+1$-th derivative is zero.

By mathematical induction it is easy to show that
\begin{equation} \label{app3}
  \begin{split}
    \frac{d^m}{dz^m}\ln\left(\frac{1+z}{1-z}\right) &= (m-1)! \left[\frac{(-1)^{m-1}}{(1+z)^m} + \frac{1}{(1-z)^m}\right] \\
     & =\frac{(m-1)!}{(1-z^2)^m} \sum_{s=0}^{m}\left[ {{m}\choose{s}}\left((-1)^{m-1+s}+ 1\right)z^s\right],
  \end{split}
\end{equation}
where for the last equality the binomial theorem
\begin{equation}
  (\alpha+\beta)^n = \sum_{s=0}^{n}{{n}\choose{s}}\alpha^{n-s}\beta^s
\end{equation}
has been employed. Another useful relation that can be found in the bibliography \cite{Heisk} is the following:
\begin{equation}
  P_\ell(z) = \frac{1}{2^\ell} \sum_{p=0}^{r} \frac{(-1)^p(2\ell-2 p)!}{p! (\ell-p)!(\ell-2p)!}z^{\ell-2p}
\end{equation}
where $r$ is the greatest integer that is smaller than $\frac{\ell}{2}$, i.e. it is $r=\frac{\ell}{2}$ if $\ell$ is even, else it is $r=\frac{\ell-1}{2}$.
Since,
\begin{equation}
  \frac{d^m}{dz^m} z^n =  z^{n-m} \prod_{j=1}^{m}(n-j+1) = \frac{(n+1-m)}{n+1} (n+2-m)_m\, z^{n-m}
\end{equation}
where in the last equality we generalize the result with the help of the Pochhammer symbol $(x)_n = \Gamma(x+n)/\Gamma(x)$. As a result, we can now write
\begin{equation} \label{app4}
  \frac{d^{\ell+1-m}}{dz^{\ell+1-m}}P_{\ell}(z) = \frac{1}{2^\ell} \sum_{p=0}^{r}\frac{(-1)^p(2\ell-2 p)!}{p! (\ell-p)!(\ell-2p)!}  \frac{(m-2 p)}{\ell-2 p+1} (m+1-2 p)_{(\ell-m+1)} z^{m-2p-1}.
\end{equation}
Due to the properties of $Q^{k}_\ell(z)$ we know that
\begin{equation} \label{app5}
  A_0(\ell) = (1-z^2)^{(\ell+1)/2} Q^{\ell+1}_{\ell}(z), \quad \ell \in \mathbb{N}
\end{equation}
is constant in respect to $z$ \cite{Susc}. By substitution of \eqref{app1}-\eqref{app3} and \eqref{app4} in \eqref{app5} we see that we have successive products of polynomials in $z$ that result into a constant value, hence for the calculation of this value we need only consider the constant terms of each (e.g. the $s=0$ of \eqref{app3} and $p=\frac{m-1}{2}$ of \eqref{app4}). A straightforward calculation yields
\begin{equation}\label{valueofA}
  \begin{split}
    A_0(\ell) = \frac{(-1)^{\ell+1}}{2^{\ell+1}} \sum_{m=1}^{\ell+1} \Bigg[ \frac{(-1)^{(m-1)/2}(\ell+1)! (2\,\ell-m+1)!(m-1)!}{m!(\ell+1-m)! \Gamma(\frac{1+m}{2})\Gamma(\ell+\frac{3-m}{2})} \left((-1)^{m-1}+1\right)\Bigg],
  \end{split}
\end{equation}
where we have also used the equality
\begin{equation}
  \frac{(2)_{(1+\ell-m)}}{\ell+2-m} = (\ell+1-m)! .
\end{equation}
Because of the $\left((-1)^{m-1}+1\right)$ multiplication term, it is clear that only the odd values of $m$ contribute in the summation in \eqref{valueofA}. Hence, the $\Gamma$ functions appearing in it can be substituted by factorials due to the fact that no half integer values are produced in the arguments. Furthermore, by setting $m=2\, n+1$ we can simplify \eqref{valueofA} and obtain
\begin{equation}\label{valueofAfinal}
  A_0(\ell) = \frac{(-1)^{\ell+1}(\ell+1)!}{2^{\ell}} \sum_{n=0}^{r} \frac{(-1)^n(2\, \ell-2\, n)!(2\, n)!}{(2\, n+1)! (\ell-2\, n)! (\ell-n)!\, n!},
\end{equation}
where now $r=\ell/2$ if $\ell$ is even, else $r=\frac{\ell-1}{2}$.

We can now write the part involving the summation in \eqref{valueofAfinal} as
\begin{equation}
S= \sum_{n=0}^{r} \frac{(-1)^n(2\, \ell-2\, n)!(2\, n)!}{(2\, n+1)! (\ell-2\, n)! (\ell-n)!\, n!} = \frac{(2 \ell)!}{(\ell!)^2} \sum^{r}_{n=0} \frac{(-1)^n \binom{\ell}{2 n}\binom{\ell}{n}}{(2n+1)\binom{2\ell}{2 n}}.
\end{equation}
Then, it can be shown, by the definition of the generalized hypergeometric function ${}_3F_2$ \cite{NIST}
\begin{equation}
  {}_3F_2(a,b,c;d,e;z) = \sum_{n=0}^{\infty} \frac{(a)_n (b)_n (c)_n}{(d)_n (e)_n} \frac{z^n}{n!}, \quad d, e \notin \mathbb{Z}_{-}\cup\{0\}
\end{equation}
and the properties of the Pochhammer symbols, that $S$ can be written as
\begin{equation}
  S = {}_{3}F_2 (\frac{1}{2},\frac{1-\ell}{2},-\frac{\ell}{2};\frac{3}{2}, \frac{1}{2}-\ell;1).
\end{equation}
If we now use the relation presented in \cite{Exton}, which is based on Saalsch\"utz' theorem,
\begin{equation}
  {}_3F_2(a,b,-n;c,a+b-c-n+1;1) = \frac{(c-a)_n (c-b)_n}{(c)_n (c-a-b)_n}
\end{equation}
together with the definition for the Pochhammer symbol $(a)_n= \Gamma(a+n)/\Gamma(n)$, the expression $\Gamma(\frac{n}{2}) = (n-2)!! \sqrt{\pi}/2^{\frac{n-1}{2}}$ for half integers $\frac{n}{2}$ and relation $(2\ell-1)!! = \frac{(2\ell)!}{2^\ell \ell!}$, then the summation results in
\begin{equation}
  S= \frac{2^{2\ell}}{\ell+1} .
\end{equation}
Finally, substitution of $S$ into \eqref{valueofAfinal} leads to
\begin{equation}\label{valueofAfinal2}
  A_0(\ell) = (-1)^{\ell+1} 2^{\ell} \ell! .
\end{equation}
The same result for $A_0(\ell)$ can be produced in a straightforward manner by using the expression for the Hobson's function $\mathcal{Q}^{k}_\ell(z)$ given in Eq. (6.17) of Ref. \cite{Szmytkowski} together with the relationships displayed in our Eqs. \eqref{relationbttwo} and \eqref{app5}.

\section*{Acknowledgements} The author acknowledges financial support by FONDECYT postdoctoral grant No. 3150016.

\end{document}